\documentclass[11pt]{article}

\usepackage[T1]{fontenc}
\usepackage[utf8]{inputenc}
\usepackage{lmodern}
\usepackage[english]{babel}
\usepackage{amsmath,amssymb,amsthm}
\usepackage{graphicx}
\usepackage{xcolor}
\usepackage{caption}
\usepackage{hyperref}
\usepackage[all]{hypcap}
\usepackage[a4paper,margin=2.8cm]{geometry}
\usepackage{tikz}
\usetikzlibrary{calc,angles,quotes,arrows.meta}
\usepackage{booktabs}
\usepackage{array}
\usepackage{microtype}
\usepackage{needspace}
\usepackage{enumitem}

\setlist[enumerate]{leftmargin=2.1em,itemsep=2pt,topsep=4pt}
\setlist[itemize]{leftmargin=2.1em,itemsep=2pt,topsep=4pt}
\hypersetup{colorlinks=true,linkcolor=blue!55!black,citecolor=blue!55!black,urlcolor=blue!55!black,pdftitle={Minimality of the double kite among strict two-pole unit-triangle configurations},pdfauthor={Mike Winkler}}

\newtheorem{theorem}{Theorem}[section]
\newtheorem{lemma}[theorem]{Lemma}
\newtheorem{proposition}[theorem]{Proposition}
\newtheorem{corollary}[theorem]{Corollary}
\theoremstyle{definition}
\newtheorem{example}[theorem]{Example}
\theoremstyle{definition}
\newtheorem{definition}[theorem]{Definition}
\theoremstyle{remark}
\newtheorem{remark}[theorem]{Remark}
\newcommand{\Area}{\operatorname{area}}
\newcommand{\Int}{\operatorname{int}}
\newcommand{\T}{\mathcal{T}}
\newcommand{\R}{\mathbb{R}}

\title{Minimality of the Double Kite among Strict Two-Pole Configurations of Unit Triangles}

\author{Mike Winkler\\[2pt]
	\small Faculty of Mathematics\\
	\small Ruhr University Bochum, Germany\\ \small \href{mailto:mike.winkler@ruhr-uni-bochum.de}{mike.winkler@ruhr-uni-bochum.de}}
	
\date{November 11, 2026}

	\long\def\DKFigureFortyTwo{%
		\begin{figure}[htbp]
			\centering
			\begin{minipage}[b]{.5\textwidth}\centering
				\begin{tikzpicture}
					[y=0.5pt, x=0.5pt, yscale=-0.9, xscale=0.9]
					\draw[line width=0.3pt]
					(460.7987,859.0503) -- (417.1412,862.1912)
					(417.1412,862.1912) -- (441.6901,898.4293)
					(441.6901,898.4293) -- (398.0326,901.5702)
					(441.6901,898.4293) -- (460.7987,859.0503)
					(460.7987,859.0503) -- (485.3476,895.2883)
					(485.3476,895.2883) -- (441.6901,898.4293)
					(485.3476,895.2883) -- (477.4744,938.3448)
					(518.6990,923.6349) -- (477.4744,938.3448)
					(437.7535,919.9575) -- (441.6901,963.5504)
					(441.6901,963.5504) -- (401.9692,945.1632)
					(401.9692,945.1632) -- (437.7535,919.9575)
					(436.2498,822.8122) -- (460.7987,859.0503)
					(485.3476,895.2883) -- (518.6990,923.6349)
					(545.6808,1006.0703) -- (521.1319,1042.3084)
					(521.1319,1042.3084) -- (502.0233,1002.9294)
					(502.0233,1002.9294) -- (477.4744,1039.1675)
					(502.0233,1002.9294) -- (545.6808,1006.0703)
					(545.6808,1006.0703) -- (526.5722,966.6914)
					(526.5722,966.6914) -- (502.0233,1002.9294)
					(526.5722,966.6914) -- (485.3476,951.9815)
					(518.6990,923.6349) -- (485.3476,951.9815)
					(481.4110,995.5745) -- (445.6267,970.3688)
					(445.6267,970.3688) -- (441.6901,1013.9618)
					(441.6901,1013.9618) -- (481.4110,995.5745)
					(564.7894,1045.4494) -- (545.6808,1006.0703)
					(526.5722,966.6914) -- (518.6990,923.6349)
					(477.4744,1039.1675) -- (521.1319,1042.3084)
					(521.1319,1042.3084) -- (564.7894,1045.4494)
					(436.2498,822.8122) -- (417.1412,862.1912)
					(417.1412,862.1912) -- (398.0326,901.5702)
					(398.0326,901.5702) -- (401.9692,945.1632)
					(401.9692,945.1632) -- (405.9058,988.7561)
					(405.9058,988.7561) -- (441.6901,1013.9618)
					(441.6901,1013.9618) -- (477.4744,1039.1675)
					(405.9058,988.7561) -- (441.6901,963.5504)
					(441.6901,963.5504) -- (477.4744,938.3448)
					(477.4744,938.3448) -- (437.7535,919.9575)
					(437.7535,919.9575) -- (398.0326,901.5702)
					(485.3476,951.9815) -- (481.4110,995.5745)
					(481.4110,995.5745) -- (477.4744,1039.1675)
					(405.9058,988.7561) -- (445.6267,970.3688)
					(445.6267,970.3688) -- (485.3476,951.9815)
					(526.5722,966.6914) -- (545.6809,1006.0703)
					(545.6808,1006.0703) -- (564.7894,1045.4494)
					(693.3291,822.8122) -- (649.6716,825.9532)
					(649.6716,825.9532) -- (606.0140,829.0941)
					(649.6716,825.9532) -- (668.7802,786.5742)
					(668.7802,786.5742) -- (625.1227,789.7151)
					(625.1227,789.7151) -- (644.2313,750.3361)
					(625.1227,789.7151) -- (649.6716,825.9532)
					(606.0140,829.0941) -- (625.1227,789.7151)
					(606.0140,829.0941) -- (572.6627,800.7475)
					(564.7895,843.8039) -- (572.6627,800.7475)
					(608.4470,775.5418) -- (568.7261,757.1545)
					(568.7261,757.1545) -- (604.5104,731.9488)
					(604.5104,731.9488) -- (608.4470,775.5418)
					(606.0140,829.0941) -- (564.7895,843.8039)
					(479.9074,825.9532) -- (460.7988,786.5742)
					(460.7988,786.5742) -- (504.4563,789.7151)
					(504.4563,789.7151) -- (485.3476,750.3361)
					(504.4563,789.7151) -- (479.9074,825.9532)
					(479.9074,825.9532) -- (523.5649,829.0941)
					(523.5649,829.0941) -- (504.4563,789.7151)
					(523.5649,829.0941) -- (556.9163,800.7475)
					(564.7895,843.8039) -- (556.9163,800.7475)
					(521.1320,775.5418) -- (560.8529,757.1545)
					(560.8529,757.1545) -- (525.0686,731.9488)
					(525.0686,731.9488) -- (521.1320,775.5418)
					(523.5649,829.0941) -- (564.7895,843.8039)
					(485.3476,750.3361) -- (460.7988,786.5742)
					(693.3291,822.8122) -- (668.7802,786.5742)
					(668.7802,786.5742) -- (644.2313,750.3361)
					(644.2313,750.3361) -- (604.5104,731.9488)
					(604.5104,731.9488) -- (564.7895,713.5615)
					(564.7895,713.5615) -- (525.0686,731.9488)
					(525.0686,731.9488) -- (485.3476,750.3361)
					(564.7895,713.5615) -- (568.7261,757.1545)
					(568.7261,757.1545) -- (572.6627,800.7475)
					(572.6627,800.7475) -- (608.4470,775.5418)
					(608.4470,775.5418) -- (644.2313,750.3361)
					(556.9163,800.7475) -- (521.1320,775.5418)
					(521.1320,775.5418) -- (485.3476,750.3361)
					(564.7895,713.5615) -- (560.8529,757.1545)
					(560.8529,757.1545) -- (556.9163,800.7475)
					(608.4470,1042.3084) -- (652.1045,1039.1675)
					(668.7802,859.0503) -- (712.4377,862.1912)
					(712.4377,862.1912) -- (687.8888,898.4293)
					(687.8888,898.4293) -- (731.5463,901.5702)
					(687.8888,898.4293) -- (668.7802,859.0503)
					(668.7802,859.0503) -- (644.2313,895.2883)
					(644.2313,895.2883) -- (687.8888,898.4293)
					(644.2313,895.2883) -- (652.1045,938.3448)
					(610.8799,923.6349) -- (652.1045,938.3448)
					(691.8254,919.9575) -- (687.8888,963.5504)
					(687.8888,963.5504) -- (727.6097,945.1632)
					(727.6097,945.1632) -- (691.8254,919.9575)
					(644.2313,895.2883) -- (610.8799,923.6349)
					(583.8981,1006.0703) -- (608.4470,1042.3084)
					(608.4470,1042.3084) -- (627.5556,1002.9294)
					(627.5556,1002.9294) -- (652.1045,1039.1675)
					(627.5556,1002.9294) -- (583.8981,1006.0703)
					(583.8981,1006.0703) -- (603.0068,966.6914)
					(603.0068,966.6914) -- (627.5556,1002.9294)
					(603.0068,966.6914) -- (644.2313,951.9815)
					(610.8799,923.6349) -- (644.2313,951.9815)
					(648.1679,995.5745) -- (683.9522,970.3688)
					(683.9522,970.3688) -- (687.8888,1013.9618)
					(687.8888,1013.9618) -- (648.1679,995.5745)
					(603.0068,966.6914) -- (610.8799,923.6349)
					(712.4377,862.1912) -- (731.5463,901.5702)
					(731.5463,901.5702) -- (727.6097,945.1632)
					(727.6097,945.1632) -- (723.6731,988.7561)
					(723.6731,988.7561) -- (687.8888,1013.9618)
					(687.8888,1013.9618) -- (652.1045,1039.1675)
					(723.6731,988.7561) -- (687.8888,963.5504)
					(687.8888,963.5504) -- (652.1045,938.3448)
					(652.1045,938.3448) -- (691.8254,919.9575)
					(691.8254,919.9575) -- (731.5463,901.5702)
					(644.2313,951.9815) -- (648.1679,995.5745)
					(648.1679,995.5745) -- (652.1045,1039.1675)
					(723.6731,988.7561) -- (683.9522,970.3688)
					(683.9522,970.3688) -- (644.2313,951.9815)
					(460.7988,786.5742) -- (436.2498,822.8122)
					(479.9074,825.9532) -- (436.2498,822.8122)
					(668.7802,859.0503) -- (693.3291,822.8122)
					(712.4377,862.1912) -- (693.3291,822.8122)
					(583.8981,1006.0703) -- (564.7894,1045.4494)
					(608.4470,1042.3084) -- (564.7894,1045.4494);
				\end{tikzpicture}
				\par\small (a) Harborth (1991)
			\end{minipage}\hfill
			\begin{minipage}[b]{.5\textwidth}\centering
				\begin{tikzpicture}
					[y=0.5pt, x=0.5pt, yscale=-0.9, xscale=0.9]
					\draw[line width=0.3pt]
					(79.7254,860.1348) -- (36.1972,864.7321)
					(36.1972,864.7321) -- (61.9427,900.1300)
					(61.9427,900.1300) -- (18.4145,904.7273)
					(61.9427,900.1300) -- (79.7254,860.1348)
					(79.7254,860.1348) -- (105.4709,895.5327)
					(105.4709,895.5327) -- (61.9427,900.1300)
					(105.4709,895.5327) -- (99.0402,938.8281)
					(139.7505,922.7496) -- (99.0402,938.8281)
					(58.7273,921.7777) -- (64.1177,965.2149)
					(64.1177,965.2149) -- (23.8048,948.1645)
					(23.8048,948.1645) -- (58.7273,921.7777)
					(53.9799,824.7369) -- (79.7254,860.1348)
					(105.4709,895.5327) -- (139.7505,922.7496)
					(169.4705,1004.2379) -- (146.1456,1041.2756)
					(146.1456,1041.2756) -- (125.7324,1002.5568)
					(125.7325,1002.5568) -- (102.4076,1039.5945)
					(125.7325,1002.5568) -- (169.4705,1004.2379)
					(169.4705,1004.2379) -- (149.0573,965.5191)
					(149.0573,965.5191) -- (125.7325,1002.5568)
					(149.0573,965.5191) -- (107.3645,952.1943)
					(139.7505,922.7496) -- (107.3645,952.1943)
					(104.8860,995.8944) -- (68.2798,971.8980)
					(68.2798,971.8980) -- (65.8014,1015.5981)
					(65.8014,1015.5981) -- (104.8860,995.8944)
					(189.8837,1042.9567) -- (169.4705,1004.2379)
					(149.0573,965.5191) -- (139.7505,922.7496)
					(102.4076,1039.5945) -- (146.1456,1041.2756)
					(146.1456,1041.2756) -- (189.8837,1042.9567)
					(53.9800,824.7369) -- (36.1972,864.7321)
					(36.1972,864.7321) -- (18.4145,904.7273)
					(18.4145,904.7273) -- (23.8048,948.1645)
					(23.8048,948.1645) -- (29.1952,991.6017)
					(29.1952,991.6017) -- (65.8014,1015.5981)
					(65.8014,1015.5981) -- (102.4076,1039.5945)
					(29.1952,991.6017) -- (64.1177,965.2149)
					(64.1177,965.2149) -- (99.0402,938.8281)
					(99.0402,938.8281) -- (58.7273,921.7777)
					(58.7273,921.7777) -- (18.4145,904.7273)
					(107.3645,952.1943) -- (104.8860,995.8944)
					(104.8860,995.8944) -- (102.4076,1039.5945)
					(29.1952,991.6017) -- (68.2798,971.8980)
					(68.2798,971.8980) -- (107.3644,952.1943)
					(300.0420,860.1348) -- (343.5703,864.7321)
					(343.5703,864.7321) -- (317.8248,900.1300)
					(317.8248,900.1300) -- (361.3530,904.7273)
					(317.8248,900.1300) -- (300.0420,860.1348)
					(300.0420,860.1348) -- (274.2965,895.5327)
					(274.2965,895.5327) -- (317.8248,900.1300)
					(274.2965,895.5327) -- (280.7273,938.8281)
					(240.0170,922.7496) -- (280.7273,938.8281)
					(321.0401,921.7777) -- (315.6498,965.2149)
					(315.6498,965.2149) -- (355.9627,948.1645)
					(355.9627,948.1645) -- (321.0402,921.7777)
					(325.7875,824.7369) -- (300.0420,860.1348)
					(274.2965,895.5327) -- (240.0170,922.7496)
					(210.2970,1004.2379) -- (233.6218,1041.2756)
					(233.6218,1041.2756) -- (254.0350,1002.5568)
					(254.0350,1002.5568) -- (277.3599,1039.5945)
					(254.0350,1002.5568) -- (210.2970,1004.2379)
					(210.2970,1004.2379) -- (230.7102,965.5191)
					(230.7102,965.5191) -- (254.0350,1002.5568)
					(230.7102,965.5191) -- (272.4030,952.1943)
					(240.0170,922.7496) -- (272.4030,952.1943)
					(274.8814,995.8944) -- (311.4876,971.8980)
					(311.4876,971.8980) -- (313.9661,1015.5981)
					(313.9661,1015.5981) -- (274.8814,995.8944)
					(230.7102,965.5191) -- (240.0170,922.7496)
					(277.3599,1039.5945) -- (233.6218,1041.2756)
					(325.7875,824.7369) -- (343.5703,864.7321)
					(343.5703,864.7321) -- (361.3530,904.7273)
					(361.3530,904.7273) -- (355.9627,948.1645)
					(355.9627,948.1645) -- (350.5723,991.6017)
					(350.5723,991.6017) -- (313.9661,1015.5981)
					(313.9661,1015.5981) -- (277.3599,1039.5945)
					(350.5723,991.6017) -- (315.6498,965.2149)
					(315.6498,965.2149) -- (280.7273,938.8281)
					(280.7273,938.8281) -- (321.0401,921.7777)
					(321.0401,921.7777) -- (361.3530,904.7273)
					(272.4030,952.1943) -- (274.8815,995.8944)
					(274.8814,995.8944) -- (277.3599,1039.5945)
					(350.5723,991.6017) -- (311.4877,971.8980)
					(311.4876,971.8980) -- (272.4030,952.1943)
					(193.1509,803.0969) -- (167.4902,767.6375)
					(167.4902,767.6375) -- (149.6118,807.5900)
					(149.6118,807.5900) -- (123.9511,772.1306)
					(149.6118,807.5900) -- (193.1509,803.0969)
					(193.1509,803.0969) -- (175.2725,843.0495)
					(175.2725,843.0495) -- (149.6118,807.5900)
					(175.2725,843.0495) -- (134.5239,859.0305)
					(168.7381,886.3293) -- (134.5239,859.0305)
					(129.2375,815.5806) -- (94.2519,841.8837)
					(94.2519,841.8837) -- (88.9655,798.4337)
					(88.9655,798.4337) -- (129.2375,815.5806)
					(230.1557,841.8837) -- (255.8164,877.3431)
					(230.1557,841.8837) -- (186.6166,846.3768)
					(186.6166,846.3768) -- (204.4950,806.4243)
					(204.4950,806.4243) -- (230.1557,841.8837)
					(204.4950,806.4243) -- (245.2436,790.4432)
					(204.4950,806.4243) -- (211.0293,763.1444)
					(250.5300,833.8932) -- (285.5156,807.5901)
					(285.5156,807.5901) -- (290.8020,851.0400)
					(290.8020,851.0400) -- (250.5300,833.8932)
					(255.8164,877.3431) -- (212.2773,881.8362)
					(212.2773,881.8362) -- (230.1557,841.8837)
					(212.2773,881.8362) -- (186.6166,846.3768)
					(175.2725,843.0495) -- (168.7382,886.3293)
					(212.2773,881.8362) -- (168.7381,886.3293)
					(211.0293,763.1444) -- (193.1509,803.0969)
					(186.6166,846.3768) -- (168.7381,886.3293)
					(285.5156,807.5901) -- (245.2436,790.4432)
					(245.2436,790.4432) -- (211.0293,763.1444)
					(211.0293,763.1444) -- (167.4902,767.6375)
					(167.4902,767.6375) -- (123.9511,772.1306)
					(123.9511,772.1306) -- (88.9655,798.4337)
					(134.5239,859.0305) -- (129.2375,815.5806)
					(129.2375,815.5806) -- (123.9511,772.1306)
					(290.8020,851.0400) -- (255.8164,877.3431)
					(255.8164,877.3431) -- (250.5300,833.8932)
					(250.5300,833.8932) -- (245.2436,790.4432)
					(285.5156,807.5901) -- (325.7875,824.7369)
					(325.7875,824.7369) -- (290.8020,851.0400)
					(134.5239,859.0305) -- (94.2519,841.8837)
					(94.2519,841.8837) -- (53.9800,824.7369)
					(53.9800,824.7369) -- (88.9655,798.4337)
					(210.2970,1004.2379) -- (189.8837,1042.9567)
					(189.8837,1042.9567) -- (233.6218,1041.2756);
				\end{tikzpicture}
				\par\small (b) Asymmetric ring (2017 catalog)
			\end{minipage}\hfill
			\caption{Two strict four-regular configurations, each containing forty-two unit triangles. (a) Harborth's symmetric ring of three double kites~\cite{Harborth1991}, also shown as $T_{42}$ in \cite[Fig.~3]{HarborthMoeller2002}. (b) An asymmetric ring of two double kites and one reverse double kite, constructed by Winkler, Dinkelacker, and Vogel in 2017 and documented in \cite[v2, Fig.~1a]{Winkler2017}. The components called \emph{double kites} and \emph{reverse double kites} are shown in \cite[v2, Fig.~9]{Winkler2017}.}
			\label{fig:fortytwo}
		\end{figure}
	}
	
	\long\def\DKKiteIntro{%
		\begin{tikzpicture}[line join=round]
			\filldraw[draw=black,line width=0.42pt,fill=orange!36] (0.00000000,0.00000000) -- (-0.25000000,0.96824584) -- (-0.96352549,0.26761657) -- cycle;
			\filldraw[draw=black,line width=0.42pt,fill=orange!36] (-0.50000000,1.93649167) -- (-0.25000000,0.96824584) -- (-1.21352549,1.23586240) -- cycle;
			\filldraw[draw=black,line width=0.42pt,fill=orange!36] (-1.92705098,0.53523313) -- (-0.96352549,0.26761657) -- (-1.21352549,1.23586240) -- cycle;
			\filldraw[draw=black,line width=0.42pt,fill=orange!36] (0.00000000,0.00000000) -- (0.25000000,0.96824584) -- (0.96352549,0.26761657) -- cycle;
			\filldraw[draw=black,line width=0.42pt,fill=orange!36] (0.50000000,1.93649167) -- (0.25000000,0.96824584) -- (1.21352549,1.23586240) -- cycle;
			\filldraw[draw=black,line width=0.42pt,fill=orange!36] (1.92705098,0.53523313) -- (0.96352549,0.26761657) -- (1.21352549,1.23586240) -- cycle;
			\filldraw[draw=black,line width=0.42pt,fill=orange!36] (-0.50000000,1.93649167) -- (0.50000000,1.93649167) -- (0.00000000,2.80251708) -- cycle;
			\fill[red!80!black] (-1.92705098,0.53523313) circle (2.2pt);
			\fill[red!80!black] (1.92705098,0.53523313) circle (2.2pt);
			\fill[red!80!black] (0.00000000,2.80251708) circle (2.2pt);
		\end{tikzpicture}
	}
	
	\long\def\DKDoubleKiteIntro{%
		\begin{tikzpicture}[line join=round]
			\filldraw[draw=black,line width=0.42pt,fill=orange!36] (0.00000000,0.00000000) -- (-0.25000000,0.96824584) -- (-0.96352549,0.26761657) -- cycle;
			\filldraw[draw=black,line width=0.42pt,fill=orange!36] (-0.50000000,1.93649167) -- (-0.25000000,0.96824584) -- (-1.21352549,1.23586240) -- cycle;
			\filldraw[draw=black,line width=0.42pt,fill=orange!36] (-1.92705098,0.53523313) -- (-0.96352549,0.26761657) -- (-1.21352549,1.23586240) -- cycle;
			\filldraw[draw=black,line width=0.42pt,fill=orange!36] (0.00000000,0.00000000) -- (0.25000000,0.96824584) -- (0.96352549,0.26761657) -- cycle;
			\filldraw[draw=black,line width=0.42pt,fill=orange!36] (0.50000000,1.93649167) -- (0.25000000,0.96824584) -- (1.21352549,1.23586240) -- cycle;
			\filldraw[draw=black,line width=0.42pt,fill=orange!36] (1.92705098,0.53523313) -- (0.96352549,0.26761657) -- (1.21352549,1.23586240) -- cycle;
			\filldraw[draw=black,line width=0.42pt,fill=orange!36] (-0.50000000,1.93649167) -- (0.50000000,1.93649167) -- (0.00000000,2.80251708) -- cycle;
			\filldraw[draw=black,line width=0.42pt,fill=blue!23] (2.76587699,2.35082420) -- (1.85058349,2.75361176) -- (2.65705449,3.34488540) -- cycle;
			\filldraw[draw=black,line width=0.42pt,fill=blue!23] (0.93528999,3.15639931) -- (1.85058349,2.75361176) -- (1.74176099,3.74767295) -- cycle;
			\filldraw[draw=black,line width=0.42pt,fill=blue!23] (2.54823199,4.33894660) -- (2.65705449,3.34488540) -- (1.74176099,3.74767295) -- cycle;
			\filldraw[draw=black,line width=0.42pt,fill=blue!23] (2.76587699,2.35082420) -- (1.76999649,2.26014875) -- (2.34646399,1.44302867) -- cycle;
			\filldraw[draw=black,line width=0.42pt,fill=blue!23] (0.77411600,2.16947330) -- (1.76999649,2.26014875) -- (1.35058349,1.35235322) -- cycle;
			\filldraw[draw=black,line width=0.42pt,fill=blue!23] (1.92705098,0.53523313) -- (2.34646399,1.44302867) -- (1.35058349,1.35235322) -- cycle;
			\filldraw[draw=black,line width=0.42pt,fill=blue!23] (0.93528999,3.15639931) -- (0.77411600,2.16947330) -- (0.00000000,2.80251708) -- cycle;
			\fill[red!80!black] (-1.92705098,0.53523313) circle (2.2pt);
			\fill[red!80!black] (2.54823199,4.33894660) circle (2.2pt);
		\end{tikzpicture}
	}
	
	\long\def\DKReverseDoubleKiteIntro{%
		\begin{tikzpicture}[line join=round]
			\filldraw[draw=black,line width=0.42pt,fill=orange!36] (0.00000000,0.00000000) -- (-0.25000000,0.96824584) -- (-0.96352549,0.26761657) -- cycle;
			\filldraw[draw=black,line width=0.42pt,fill=orange!36] (-0.50000000,1.93649167) -- (-0.25000000,0.96824584) -- (-1.21352549,1.23586240) -- cycle;
			\filldraw[draw=black,line width=0.42pt,fill=orange!36] (-1.92705098,0.53523313) -- (-0.96352549,0.26761657) -- (-1.21352549,1.23586240) -- cycle;
			\filldraw[draw=black,line width=0.42pt,fill=orange!36] (0.00000000,0.00000000) -- (0.25000000,0.96824584) -- (0.96352549,0.26761657) -- cycle;
			\filldraw[draw=black,line width=0.42pt,fill=orange!36] (0.50000000,1.93649167) -- (0.25000000,0.96824584) -- (1.21352549,1.23586240) -- cycle;
			\filldraw[draw=black,line width=0.42pt,fill=orange!36] (1.92705098,0.53523313) -- (0.96352549,0.26761657) -- (1.21352549,1.23586240) -- cycle;
			\filldraw[draw=black,line width=0.42pt,fill=orange!36] (-0.50000000,1.93649167) -- (0.50000000,1.93649167) -- (0.00000000,2.80251708) -- cycle;
			\filldraw[draw=black,line width=0.42pt,fill=blue!23] (1.92705098,3.33775021) -- (2.17705098,2.36950437) -- (2.89057647,3.07013364) -- cycle;
			\filldraw[draw=black,line width=0.42pt,fill=blue!23] (2.42705098,1.40125854) -- (2.17705098,2.36950437) -- (3.14057647,2.10188781) -- cycle;
			\filldraw[draw=black,line width=0.42pt,fill=blue!23] (3.85410197,2.80251708) -- (2.89057647,3.07013364) -- (3.14057647,2.10188781) -- cycle;
			\filldraw[draw=black,line width=0.42pt,fill=blue!23] (1.92705098,3.33775021) -- (1.67705098,2.36950437) -- (0.96352549,3.07013364) -- cycle;
			\filldraw[draw=black,line width=0.42pt,fill=blue!23] (1.42705098,1.40125854) -- (1.67705098,2.36950437) -- (0.71352549,2.10188781) -- cycle;
			\filldraw[draw=black,line width=0.42pt,fill=blue!23] (0.00000000,2.80251708) -- (0.96352549,3.07013364) -- (0.71352549,2.10188781) -- cycle;
			\filldraw[draw=black,line width=0.42pt,fill=blue!23] (2.42705098,1.40125854) -- (1.42705098,1.40125854) -- (1.92705098,0.53523313) -- cycle;
			\fill[red!80!black] (-1.92705098,0.53523313) circle (2.2pt);
			\fill[red!80!black] (3.85410197,2.80251708) circle (2.2pt);
		\end{tikzpicture}
	}
	
	\long\def\DKDoubleKite{%
		\begin{tikzpicture}[scale=1,line join=round]
			\filldraw[draw=black,line width=0.35pt,fill=orange!36] (0.00000000,0.00000000) -- (-0.25000000,0.96824584) -- (-0.96352549,0.26761657) -- cycle;
			\filldraw[draw=black,line width=0.35pt,fill=orange!36] (-0.50000000,1.93649167) -- (-0.25000000,0.96824584) -- (-1.21352549,1.23586240) -- cycle;
			\filldraw[draw=black,line width=0.35pt,fill=orange!36] (-1.92705098,0.53523313) -- (-0.96352549,0.26761657) -- (-1.21352549,1.23586240) -- cycle;
			\filldraw[draw=black,line width=0.35pt,fill=orange!36] (0.00000000,0.00000000) -- (0.25000000,0.96824584) -- (0.96352549,0.26761657) -- cycle;
			\filldraw[draw=black,line width=0.35pt,fill=orange!36] (0.50000000,1.93649167) -- (0.25000000,0.96824584) -- (1.21352549,1.23586240) -- cycle;
			\filldraw[draw=black,line width=0.35pt,fill=orange!36] (1.92705098,0.53523313) -- (0.96352549,0.26761657) -- (1.21352549,1.23586240) -- cycle;
			\filldraw[draw=black,line width=0.35pt,fill=orange!36] (-0.50000000,1.93649167) -- (0.50000000,1.93649167) -- (0.00000000,2.80251708) -- cycle;
			\filldraw[draw=black,line width=0.35pt,fill=blue!23] (2.76587699,2.35082420) -- (1.85058349,2.75361176) -- (2.65705449,3.34488540) -- cycle;
			\filldraw[draw=black,line width=0.35pt,fill=blue!23] (0.93528999,3.15639931) -- (1.85058349,2.75361176) -- (1.74176099,3.74767295) -- cycle;
			\filldraw[draw=black,line width=0.35pt,fill=blue!23] (2.54823199,4.33894660) -- (2.65705449,3.34488540) -- (1.74176099,3.74767295) -- cycle;
			\filldraw[draw=black,line width=0.35pt,fill=blue!23] (2.76587699,2.35082420) -- (1.76999649,2.26014875) -- (2.34646399,1.44302867) -- cycle;
			\filldraw[draw=black,line width=0.35pt,fill=blue!23] (0.77411600,2.16947330) -- (1.76999649,2.26014875) -- (1.35058349,1.35235322) -- cycle;
			\filldraw[draw=black,line width=0.35pt,fill=blue!23] (1.92705098,0.53523313) -- (2.34646399,1.44302867) -- (1.35058349,1.35235322) -- cycle;
			\filldraw[draw=black,line width=0.35pt,fill=blue!23] (0.93528999,3.15639931) -- (0.77411600,2.16947330) -- (0.00000000,2.80251708) -- cycle;
			\draw[densely dashed,gray!75,line width=.5pt] (0.00000000,2.80251708) -- (1.92705098,0.53523313);
			\fill[red!85!black] (-1.92705098,0.53523313) circle (1.6pt);
			\fill[red!85!black] (2.54823199,4.33894660) circle (1.6pt);
			\node[red!75!black,anchor=east] at (-1.92705098,0.53523313) {$b_1$};
			\node[red!75!black,anchor=west] at (2.54823199,4.33894660) {$\sigma(b_1)$};
		\end{tikzpicture}
	}
	
\begin{document}
		\maketitle
		\begin{abstract}
			A strict unit-triangle configuration is a finite connected family of closed equilateral unit triangles in which distinct triangles are disjoint or meet at a common vertex, and no vertex belongs to more than two triangles. A vertex belonging to exactly one triangle is a pole. We prove that every strict configuration with one or two poles contains at least fourteen triangles, a sharp bound for two poles attained by the double kite and the reverse double kite. A strict four-regular configuration whose contact graph has an edge cut of size at most two contains at least twenty-eight triangles. A decomposition into three strict two-pole components requires forty-two triangles; both known ring constructions attain the bound. The proof reduces minimal counterexamples to two-connected contact graphs, applies area bounds to complementary unit-sided polygons, and excludes the remaining cases by exact angle certificates. The supplement contains 13,668 stored marked records, 233 exclusion certificates, a standard-library verifier, and independent checks of the complete enumeration.
		\end{abstract}

		\section{Introduction}
		
		Harborth~\cite{Harborth1991} constructed the symmetric vertex-to-vertex arrangement of forty-two unit triangles shown in Figure~\ref{fig:fortytwo}(a). Harborth and M\"oller later described the seven- and fourteen-triangle components as $T_7$ and $T_{14}$, respectively~\cite[Fig.~1]{HarborthMoeller2002}, and exhibited the complete packing $T_{42}$ assembled from three copies of $T_{14}$~\cite[Fig.~3]{HarborthMoeller2002}. They explicitly stated the conjecture that forty-two is the smallest number of triangles in a complete vertex-to-vertex packing~\cite[p.~115]{HarborthMoeller2002}. A second, asymmetric forty-two-triangle arrangement, constructed by Winkler, Dinkelacker, and Vogel in 2017, is documented in Winkler's catalog~\cite[v2, Fig.~1a]{Winkler2017}; the same catalog reproduces the symmetric ring in Figure~1b and shows the fourteen-triangle two-pole components in Figure~9. The asymmetric arrangement does not appear in Harborth and M\"oller's 2002 paper. Both arrangements are strict: distinct triangles meet only at common vertices, each of which belongs to exactly two triangles. Each arrangement consists of six seven-triangle kites grouped into three fourteen-triangle two-pole components.
		
		\DKFigureFortyTwo
		
		Thus strict four-regular configurations exist with forty-two unit triangles, while their minimum number remains unknown.
		
\begin{figure}[htbp]
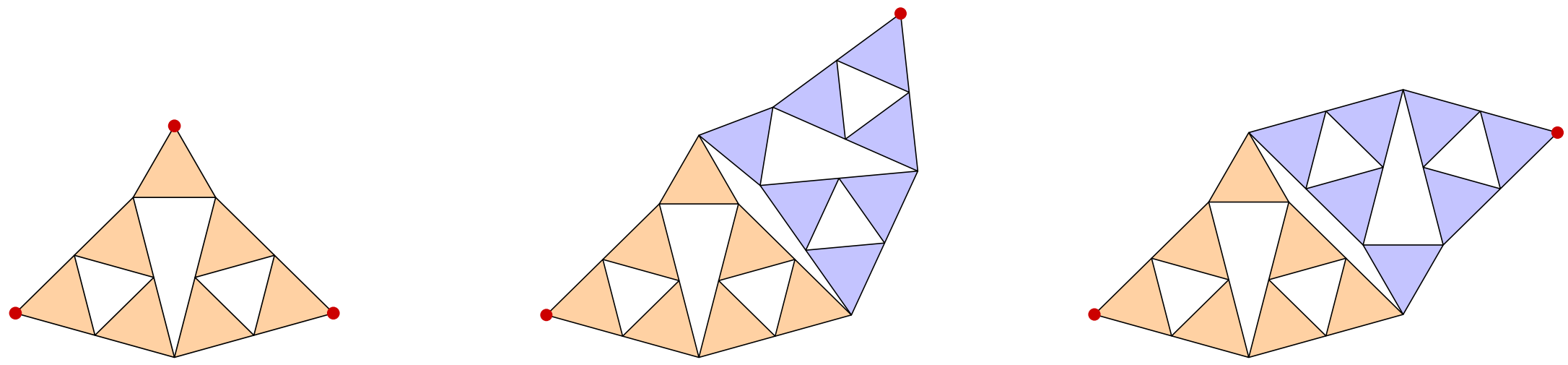

			\centering
			\begin{minipage}[b]{.29\textwidth}\centering
				\resizebox{.70\linewidth}{!}{\DKKiteIntro}
				\par\footnotesize (a) Kite: 7 triangles
			\end{minipage}\hfill
			\begin{minipage}[b]{.33\textwidth}\centering
				\resizebox{.70\linewidth}{!}{\DKDoubleKiteIntro}
				\par\footnotesize (b) Double kite: 14 triangles
			\end{minipage}\hfill
			\begin{minipage}[b]{.34\textwidth}\centering
				\resizebox{.85\linewidth}{!}{\DKReverseDoubleKiteIntro}
				\par\footnotesize (c) Reverse double kite: 14 triangles
			\end{minipage}
			\caption{A seven-triangle kite and two fourteen-triangle configurations obtained by joining two copies. Orange and blue distinguish the two constituent kites; red dots mark the poles. The drawings use the exact constructions given in Section~\ref{sec:construction}.}
			\label{fig:kite-pair}
		\end{figure}

		The configuration denoted by $T_7$ in \cite{HarborthMoeller2002} is the \emph{kite}: it consists of seven unit triangles with three free corners, each belonging to one triangle; all other corners belong to two. If edges are drawn along the triangle sides, their degrees are two and four, respectively. The names \emph{kite}, \emph{double kite}, and \emph{reverse double kite} were introduced by Winkler, Dinkelacker, and Vogel in \cite[Sec.~2 and Fig.~1]{WDV2016}. The double kite corresponds to $T_{14}$ in~\cite{HarborthMoeller2002}. Joining two kites as in Figure~\ref{fig:kite-pair} produces either double-kite configuration. Both have fourteen unit triangles, twenty-two distinct vertices, and exactly two vertices of degree two, called \emph{poles}.
		
		The examples motivate a smaller extremal question: \emph{Can a connected strict two-pole configuration contain fewer than fourteen unit triangles?} Theorem~\ref{thm:main} answers this in the negative. Two-pole configurations arise naturally by separating a strict four-regular configuration along a two-edge cut of its contact graph. Our theorem therefore also bounds four-regular configurations whose contact graphs are not three-edge-connected (Corollary~\ref{cor:edgecut}). For a cyclic union of three strict two-pole components it gives the sharp bound of forty-two triangles (Corollary~\ref{cor:ring}). The unrestricted minimum for four-regular configurations remains open.
		
		Planarity of a contact graph alone does not ensure a realization by unit triangles. Our proof reduces a minimal counterexample to a two-connected contact graph, derives necessary area inequalities for its complementary polygons, and excludes the remaining plane maps using exact angle certificates. The finite search ends at twelve triangles. We do not classify the fourteen-triangle configurations.
		
		\begin{definition}\label{def:model}
			A \emph{strict unit-triangle configuration} is a finite family $\T$ of closed equilateral triangles of side length one in $\R^2$, satisfying the following conditions:
			\begin{enumerate}
				\item For distinct $A,B\in\T$, either $A\cap B=\varnothing$, or $A\cap B$ consists of one point that is a vertex of both triangles.
				\item No point is a vertex of more than two triangles of $\T$.
				\item The union $U=\bigcup_{A\in\T}A$ is connected.
			\end{enumerate}
			A vertex occurring in exactly one triangle is a \emph{pole}. A configuration with $r$ poles is an $r$-pole configuration. Its \emph{side graph} $G(\T)$ has the distinct triangle corners as vertices and precisely the triangle sides as edges. At a pole its degree is two; at every other vertex its degree is four. A configuration without poles is called \emph{four-regular}, referring to this side graph, not to the graph of all pairs of corners at unit distance. Such a configuration is a \emph{complete vertex-to-vertex packing} in the terminology of Harborth and M\"oller~\cite{HarborthMoeller2002}.
		\end{definition}
		
		The condition of strictness excludes shared edges, intersections in the relative interior of a triangle side, and overlapping triangles. We call the triangles \emph{black} and the connected components of $\R^2\setminus U$ \emph{white faces}. Contact angles are measured in degrees; trigonometric expressions for polygon areas use radians. The marked combinatorial cases used in the proof are defined in Section~\ref{sec:census}.
		
		\begin{theorem}\label{thm:main}
			Every strict configuration with one or two poles contains at least fourteen unit triangles. For two poles this bound is attained.
		\end{theorem}
		
		Sections~\ref{sec:reduction} and~\ref{sec:polygons} develop the reduction and necessary geometric conditions. Sections~\ref{sec:census} and~\ref{sec:certs} establish the finite exclusion; Section~\ref{sec:construction} gives exact constructions attaining the bound. The cases with fourteen triangles are not enumerated.

		\section{Contact graphs and minimal-counterexample reduction}\label{sec:reduction}
		
		The \emph{contact graph} $H$ of $\T$ has one vertex for each black triangle and an edge for each common vertex of two triangles. It is connected if and only if the union of the triangles is connected. Its geometric embedding is obtained by joining the center of each triangle to its contact vertices. The two segments meeting at each contact give an embedded edge of $H$. They lie in triangle interiors except at their common endpoint and consequently cannot cross.
		
		\begin{figure}[htbp]
			\centering
			\begin{tikzpicture}[scale=1.5]
				\coordinate (O) at (0,0);
				\coordinate (A) at (1.15,0);
				\coordinate (B) at (0.575,0.99593);
				\coordinate (C) at (180:0);
				\coordinate (D) at (160:1.15);
				\coordinate (E) at (220:1.15);
				\filldraw[fill=blue!16,draw=black,line width=.7pt] (O)--(A)--(B)--cycle;
				\filldraw[fill=orange!25,draw=black,line width=.7pt] (O)--(D)--(E)--cycle;
				\fill (O) circle (1.4pt);
				\draw[thin,blue!65!black] (0.38,0) arc[start angle=0,end angle=60,radius=.38];
				\draw[thin,red!70!black] (60:.53) arc[start angle=60,end angle=160,radius=.53];
				\draw[thin,blue!65!black] (160:.36) arc[start angle=160,end angle=220,radius=.36];
				\draw[thin,red!70!black] (220:.53) arc[start angle=220,end angle=360,radius=.53];
				\node[font=\small,red!65!black] at (110:.90) {$x$};
				\node[font=\small,red!65!black] at (290:.87) {$240^\circ-x$};
				\node[font=\scriptsize] at (0.50,.32) {$60^\circ$};
				\node[font=\scriptsize] at (-.53,-.12) {$60^\circ$};
				\node[font=\small,anchor=west] at (1.4,.38) {black unit triangles};
				\draw[-{Latex[length=2mm]},gray] (1.42,.29)--(.67,.28);
			\end{tikzpicture}
			\caption{Two black triangles meeting at one contact. The white sectors satisfy $0^\circ<x<240^\circ$, and their sum is $240^\circ$. The strict inequalities follow from the prohibition of overlapping unit sides.}\label{fig:contact}
		\end{figure}
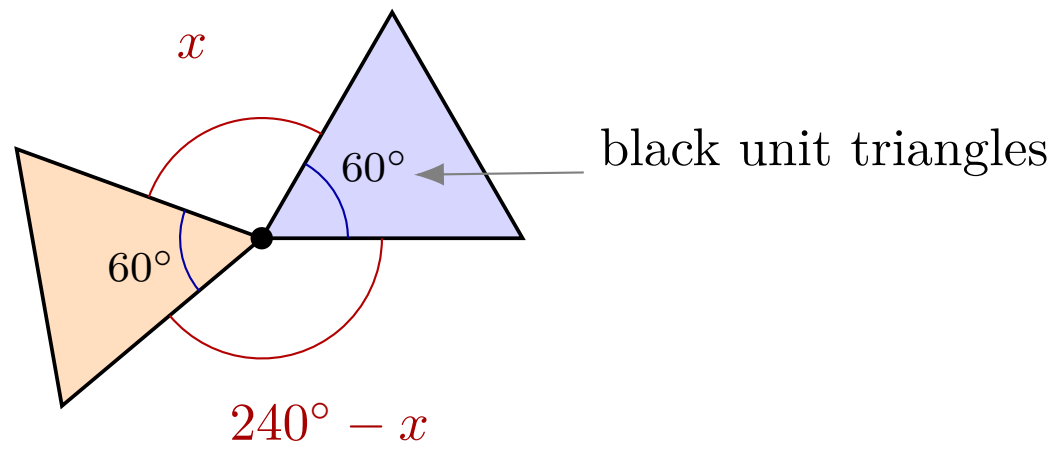
		
		\begin{lemma}[Incidence and parity]\label{lem:count}
			If $\T$ contains $t$ triangles, $v$ distinct vertices, and $r$ poles, then
			\begin{equation}\label{eq:incidence}
				3t=2v-r,\qquad r\equiv t\pmod 2.
			\end{equation}
			If $r=2$, the contact graph has $t$ vertices and $(3t-2)/2$ edges.
		\end{lemma}
		\begin{proof}
			Count the $3t$ pairs consisting of a black triangle and one of its vertices. The $v-r$ non-poles each occur twice and the $r$ poles once, so $3t=2(v-r)+r=2v-r$. Each non-pole determines exactly one contact edge, and two triangles have at most one common point. Hence non-poles and contact-graph edges are in bijection.
		\end{proof}
		
		For the lower bound it is useful to exclude one-pole configurations simultaneously. The following reduction avoids any assumption about the location of the poles in the unbounded face.
		
		\begin{lemma}[Minimal counterexample]\label{lem:minimal}
			Suppose there exists a strict configuration with one or two poles and at most thirteen triangles. Choose one with the smallest number of triangles. Then it has the following properties:
			\begin{enumerate}
				\item There are exactly two poles, lying in distinct nonadjacent triangles.
				\item Its contact graph is simple and has no bridge. In particular, it is two-connected.
				\item Its contact graph has two vertices of degree two, all other vertices have degree three, and $t\in\{4,6,8,10,12\}$.
			\end{enumerate}
		\end{lemma}
		\begin{proof}
			Each connected component obtained by deleting triangles remains a strict configuration. Its poles are the retained original poles and the contact vertices freed by the deletions.
			
			Suppose first that the configuration has a single pole in a triangle $P$. The other two vertices of $P$ meet two distinct triangles $A$ and $B$. After removing $P$, every component contains $A$ or $B$; it has one or two poles, according as it contains one or both. This contradicts minimality.
			
			If both poles lie in the same triangle $P$, that triangle is a leaf of the contact graph. Deleting it leaves a connected configuration with exactly one pole. If the pole triangles $P_1,P_2$ are adjacent, each has one further contact. Removing both triangles yields one or two nonempty components, each with one or two poles. The new poles are distinct: if their contacts with the same third triangle coincided, that point would belong to three black triangles. If they lie in different third triangles, they are distinct because each vertex belongs to at most two triangles. Both possibilities contradict minimality.
			
			Finally, suppose that an edge of the contact graph is a bridge. Cutting at its contact produces two connected strict configurations. Each acquires one new pole; the original two poles are divided between them. At least one side consequently has one or two poles, again contradicting minimality. Thus the graph has no bridge. A connected graph of maximum degree three without bridges has no cut vertex: otherwise one component at a cut vertex would be attached by only one edge. Hence the graph is two-connected. The degree sequence and the parity assertion follow from Lemma~\ref{lem:count}; the case $t=2$ is impossible for a simple contact graph.
		\end{proof}
		
		\begin{lemma}[White faces and their boundaries]\label{lem:faces}
			In the setting of Lemma~\ref{lem:minimal}, the faces of the embedded contact graph correspond to the connected components of $\R^2\setminus U$. Each such \emph{white face} has a simple closed polygonal boundary consisting of unit sides of black triangles. Let $\ell_W$ be the length of its boundary cycle in the contact graph, and let $p_W$ be the number of poles belonging to its polygonal boundary. Then this boundary has
			\begin{equation}\label{eq:face-length}
				k_W=\ell_W+p_W
			\end{equation}
			unit sides. The unbounded white face has a simple outer boundary, denoted by $Q$.
		\end{lemma}
		\begin{proof}
			Join the center of every black triangle to its contact vertices by straight spokes. The spokes inside any one triangle do not cross, and spokes from distinct triangles meet only at a prescribed contact; they give a planar embedding of the contact graph. The complement of this embedded graph consists of white regions together with the open sectors cut out by the spokes inside the black triangles. Each such black sector is adjacent to exactly one white component along its black side, or along the two consecutive black sides meeting at a pole. Attaching or removing these sectors does not join distinct white components. Thus the embedded graph faces correspond bijectively to the white faces.
			
			The graph is two-connected by Lemma~\ref{lem:minimal}, so the facial walks are cycles with pairwise distinct graph vertices (see, e.g.,~\cite{Diestel}). Different edges of a facial cycle represent different contact points: otherwise at least three distinct black triangles would contain the same point, contrary to Definition~\ref{def:model}. A facial cycle $(u_1,\ldots,u_\ell)$ determines, in each triangle $u_i$, the side joining its contacts with $u_{i-1}$ and $u_{i+1}$; if its pole lies in this face, the boundary instead traverses the two sides through that pole. At a degree-two triangle, the two center-to-contact spokes separate a sector containing its free corner and the two sides meeting there from a sector containing its remaining side. The two incident graph faces are distinct by two-connectivity; hence the pole belongs to exactly one of their white boundaries. No two nonconsecutive boundary sides can meet. Indeed, a nonconsecutive contact with a third triangle is at the unused third vertex of a degree-three triangle, not on the selected side; degree-two triangles have no third contact. Strictness forbids any other intersection of black sides. Free corners are distinct and cannot coincide with contacts. Hence the resulting boundary is a simple polygon.
			
			A triangle without a pole contributes one unit side; a triangle with its pole on the face contributes two. Therefore its polygon has $\ell_W+p_W$ unit sides. Consecutive collinear unit sides are retained as separate sides, and their common angle is $180^\circ$.
		\end{proof}

		\section{Geometric constraints on the white polygons}\label{sec:polygons}
		
		Set $a_0=\sqrt3/4$, the area of a unit equilateral triangle. For $k\ge3$, put
		\begin{equation}\label{eq:area-max}
			M(k)=\frac{k}{4}\cot\frac{\pi}{k}.
		\end{equation}
		The polygonal isoperimetric inequality states that a regular $k$-gon maximizes the area among all simple $k$-gons with fixed perimeter; in particular, a simple $k$-gon with unit sides has area at most $M(k)$ \cite{Blasjo2005}. The relevant lower bound originates in the answer by B\"or\"oczky, Kert\'esz, and Makai \cite[Theorem~1]{BKM2000} to a question of Harborth: a simple polygon with an odd number of unit sides has area at least $a_0$. Their characterization of the infimum uses a triangle whose side lengths arise from signed integral sums of the prescribed side lengths. For an odd number of unit sides this yields a nondegenerate triangle with integer sides and odd perimeter; Heron's formula gives area at least $\sqrt3/4$. The non-strict bound also applies when consecutive unit sides are collinear.
		
		\begin{lemma}[Weak area inequality]\label{lem:area}
			Let $\T$ have the structural properties (1)--(3) of Lemma~\ref{lem:minimal}. Let $k$ denote the number of unit sides of its exterior boundary $Q$, and let $f_{\rm odd}$ be the number of bounded white faces with an odd number of sides. Then
			\begin{equation}\label{eq:area-bound}
				(t+f_{\rm odd})a_0\le \Area(\Int Q)\le M(k).
			\end{equation}
		\end{lemma}
		\begin{proof}
			The interior of $Q$ is, up to sets of area zero, the disjoint union of the $t$ black triangles and the bounded white faces. Each black triangle has area $a_0$, while each bounded odd-sided white face has area at least $a_0$. This proves the first inequality. The polygonal isoperimetric inequality applied directly to the simple unit-sided polygon $Q$ gives the second inequality.
		\end{proof}
		
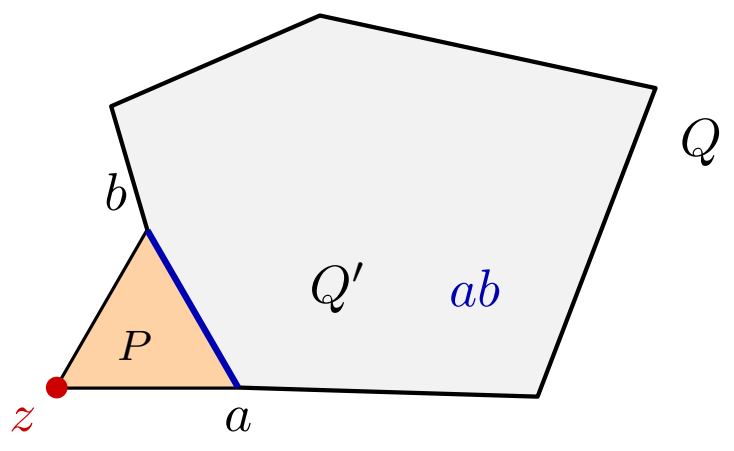
\begin{figure}[htbp]
			\centering
			\begin{tikzpicture}[scale=1.2,line join=round]
				\coordinate (P) at (0,0);
				\coordinate (A) at (1,0);
				\coordinate (B) at (.5,.8660254);
				\coordinate (C) at (2.65,-.05);
				\coordinate (D) at (3.3,1.65);
				\coordinate (E) at (1.45,2.05);
				\coordinate (F) at (.3,1.55);
				\filldraw[fill=gray!10,draw=black,thick] (P)--(A)--(C)--(D)--(E)--(F)--(B)--cycle;
				\filldraw[fill=orange!35,draw=black,thin] (P)--(A)--(B)--cycle;
				\draw[very thick,blue!70!black] (A)--(B);
				\fill[red!80!black] (P) circle (1.7pt);
				\node[below left,font=\small,red!80!black] at (P) {$z$};
				\node[below,font=\small] at (A) {$a$};
				\node[above left,font=\small] at (B) {$b$};
				\node[font=\small] at (1.55,.55) {$Q'$};
				\node[font=\small,blue!70!black,anchor=west] at (2.05,.55) {$ab$};
				\node[font=\scriptsize] at (.43,.23) {$P$};
				\node[font=\small] at (3.55,1.35) {$Q$};
			\end{tikzpicture}
			\caption{Cutting off an exterior pole ear. Replacing the boundary path $a\!\! -\!\! z\!\! -\!\! b$ by the third unit side $ab$ removes exactly the black triangle $P$, decreases the boundary length by one, and decreases the enclosed area by $a_0$. The rest of the boundary may be nonconvex.}\label{fig:ear}
		\end{figure}

		An exterior pole permits a stronger area bound, by the boundary operation shown in Figure~\ref{fig:ear}.
		
		\begin{lemma}[Exterior-pole ear inequality]\label{lem:ear}
			Let $p\in\{0,1,2\}$ be the number of poles lying on the exterior boundary $Q$ of a configuration as in Lemma~\ref{lem:minimal}. Then
			\begin{equation}\label{eq:ear-bound}
				(t-p+f_{\rm odd})a_0\le M(k-p).
			\end{equation}
		\end{lemma}
		\begin{proof}
			Consider a pole $z$ of the black triangle $P=\triangle abz$ on $Q$. The boundary follows the two consecutive sides $az$ and $zb$. The third side $ab$ is not a side of $Q$: it borders the other white face incident with $P$, which is distinct by two-connectivity. Strictness prevents any other boundary side from meeting its relative interior. Therefore $ab$ is an interior diagonal of $Q$. Replacing $a-z-b$ by $ab$ produces another simple polygon, with one fewer unit side and with area smaller by exactly $a_0$, as in Figure~\ref{fig:ear}.
			
			If both poles are exterior, their black triangles are disjoint by Lemma~\ref{lem:minimal}. Thus the second ear remains available after the first is removed; the resulting polygon is still simple. After $p$ removals the polygon encloses the other $t-p$ black triangles and all bounded white faces. Moreover, $k-p=\ell_Q\ge3$, since the exterior facial boundary of the simple two-connected contact graph is a cycle. Thus $M(k-p)$ is defined, and the lower area bound and the isoperimetric upper bound give \eqref{eq:ear-bound}.
		\end{proof}
		
		For the angular model we fix a counterclockwise cyclic order of neighbors at each vertex. A directed contact $(u,v)$, or \emph{dart}, belongs to the face obtained by following the successor of $u$ in the cyclic order at $v$. Denote by $x_{uv}$ the white angle at the contact $uv$ in that face. Angles are in degrees.
		
		\begin{lemma}[Pole elimination and reduced face angles]\label{lem:pole-elimination}
			Let $W$ be a white face with facial cycle of length $\ell$. For each pole $z$ lying on its boundary, let $P=\triangle azb$ be the corresponding black triangle, where $a,b$ are its contacts. Replacing the boundary path $a-z-b$ by the unit segment $ab$, for every such pole, produces a simple polygon $W^*$ with $\ell$ unit sides and with precisely the contact points as vertices.
			
			For a dart $(u,v)$ of $W$, let $\varepsilon_{uv}$ be the number of eliminated pole corners adjacent to its contact, and let $\alpha^*_{uv}$ denote the interior angle of $W^*$ at that point. Then $\varepsilon_{uv}\in\{0,1\}$ and
			\begin{equation}\label{eq:reduced-angles}
				\alpha^*_{uv}=\begin{cases}
					x_{uv}+60\varepsilon_{uv},&W\text{ bounded},\\
					360-x_{uv}-60\varepsilon_{uv},&W\text{ unbounded}.
				\end{cases}
			\end{equation}
			Moreover, $\sum_{(u,v)\in\partial W}\varepsilon_{uv}=2p_W$.
		\end{lemma}
		\begin{proof}
			The contact graph is two-connected. The side $ab$ borders the other face incident with the degree-two vertex $P$, hence is not part of the boundary of $W$. By strictness its relative interior does not meet any other boundary side. In a bounded white face, adjoining the interior of $P$ fills in a triangular notch; for the exterior face, replacing the two boundary sides cuts the triangular ear off the bounded region enclosed by its outer boundary. The Jordan curve theorem implies in either case that the new polygon is simple.
			
			The two pole triangles of a minimal counterexample are nonadjacent, hence disjoint. Consequently these operations do not interfere even when both poles belong to the same face. Removing one pole adds $60^\circ$ to each of the two interior angles at its contacts in the bounded-face case. For the exterior polygon it decreases the corresponding polygon angles by $60^\circ$. At a contact there can be at most one adjacent pole triangle, because the pole triangles are nonadjacent. This proves \eqref{eq:reduced-angles}; counting the two endpoints of every removed ear gives the last identity.
		\end{proof}
		
		\begin{lemma}[Contact and face equations]\label{lem:angles}
			For every edge $uv$ of the contact graph,
			\begin{equation}\label{eq:contact}
				0<x_{uv}<240,\qquad x_{uv}+x_{vu}=240.
			\end{equation}
			If a white face has a contact-graph boundary cycle of length $\ell$ and contains $p_W$ poles on its polygonal boundary, then
			\begin{equation}\label{eq:face-angle}
				\sum_{(u,v)\in\partial W}x_{uv}=
				\begin{cases}
					180(\ell-2)-120p_W,&W\text{ bounded},\\
					180(\ell+2)-120p_W,&W\text{ unbounded}.
				\end{cases}
			\end{equation}
		\end{lemma}
		\begin{proof}
			At a contact the two black angles are $60^\circ$, so the two complementary white sectors sum to $240^\circ$ (Figure~\ref{fig:contact}). A sector of zero angle would force two unit sides of different triangles to overlap, contrary to strictness. Apply Lemma~\ref{lem:pole-elimination}. The angles of the reduced $\ell$-gon sum to $180(\ell-2)$. Substituting \eqref{eq:reduced-angles} and $\sum\varepsilon_{uv}=2p_W$ gives both instances of \eqref{eq:face-angle}.
		\end{proof}
		
\begin{figure}[htbp]
			\centering
			\begin{tikzpicture}[scale=1.15,line join=round]
				\begin{scope}[shift={(0,0)}]
					\coordinate (A) at (0,0);\coordinate (B) at (1,0);\coordinate (C) at (.5,.8660254);
					\filldraw[fill=blue!20,draw=black] (A)--(B)--(C)--cycle;
					\node[font=\small] at (.5,-.35) {$B_1$};
				\end{scope}
				\begin{scope}[shift={(3.1,-.5)}]
					\coordinate (A) at (0,0);\coordinate (B) at (2,0);\coordinate (C) at (1,1.7320508);
					\coordinate (D) at (1,0);\coordinate (E) at (1.5,.8660254);\coordinate (F) at (.5,.8660254);
					\filldraw[fill=orange!32,draw=black] (A)--(D)--(F)--cycle;
					\filldraw[fill=orange!32,draw=black] (B)--(D)--(E)--cycle;
					\filldraw[fill=orange!32,draw=black] (C)--(E)--(F)--cycle;
					\draw[dashed,gray] (D)--(E)--(F)--cycle;
					\node[font=\small] at (1,-.35) {$B_3$};
					\node[font=\tiny] at (1,.59) {white};
				\end{scope}
			\end{tikzpicture}
			\caption{A single black unit triangle $B_1$ and a three-triangle block $B_3$ surrounding a white unit triangle. The central white triangle corresponds to a three-cycle of the contact graph.}\label{fig:blocks}
		\end{figure}
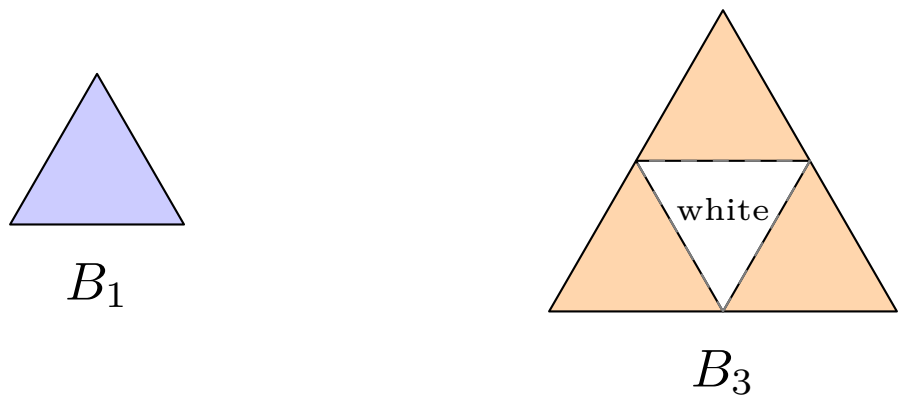

		\begin{lemma}[Three-cycles]\label{lem:three}
			Every three-cycle of the contact graph surrounds a white unit equilateral triangle. Two three-cycles cannot share a vertex.
		\end{lemma}
		\begin{proof}
			If three black triangles are pairwise in contact, their three distinct contact points are vertices of each participating black triangle in pairs. The segments joining these points are therefore three black unit sides and form an equilateral triangle. Its interior is white: a black triangle occupying it would share a full side with one of the three given triangles, while a triangle crossing its boundary would violate strictness. Hence the three white angles are $60^\circ$. The resulting arrangement is the three-triangle block shown in Figure~\ref{fig:blocks}.
			
			If two three-cycles shared a vertex in a graph of maximum degree three, they would share an edge. The two white faces on the two sides of that contact would both have angle $60^\circ$, contradicting \eqref{eq:contact}.
		\end{proof}
		
		After pole elimination, a white face is represented by the polygon $W^*$ of Lemma~\ref{lem:pole-elimination}. We may contract a corner whose interior angle is proved to be $180^\circ$, replacing its two collinear unit-sided chains by a single side with length equal to their sum. All remaining side lengths are positive integers.
		
		\begin{lemma}[Contracted-polygon conditions]\label{lem:polygon}
			Consider a simple white-face boundary after eliminating its pole corners and deleting any vertices known to have interior angle $180^\circ$. Label its vertices cyclically $v_1,\ldots,v_m$, with $L_j$ the positive integer length of the side $v_jv_{j+1}$ (indices modulo $m$).
			\begin{enumerate}
				\item It has $m\ge3$ and $2\max_j L_j<\sum_jL_j$.
				\item If $m=3$, its angles are determined by the cosine rule.
				\item If $m=4$, $L_1=L_2$, and $L_3=L_4$, then $v_2$ and $v_4$ are the apex vertices, $v_1$ and $v_3$ are the base vertices, and the two base angles are equal. This includes rhombi and nonconvex kite quadrilaterals.
				\item Suppose an angle $\alpha$ has adjacent side lengths $a,b$, and the opposite chain has lengths $r_1,\ldots,r_s$. Put $S=\sum_i r_i$ and $L_0=\max(0,2\max_i r_i-S)$. If $L_0>|a-b|$, then
				\begin{equation}\label{eq:chord}
					\arccos\frac{a^2+b^2-L_0^2}{2ab}
					\le \alpha \le
					360^\circ-\arccos\frac{a^2+b^2-L_0^2}{2ab}.
				\end{equation}
			\end{enumerate}
		\end{lemma}
		\begin{proof}
			The strict polygon inequality in (1) follows by viewing a side as one part of a noncollinear closed chain. Part (2) is the cosine rule. For (3), reflection in the line through the apex vertices exchanges the base vertices and their angles. The chain joining the neighbors of the angle in (4) has endpoint distance at least $L_0$, by the reverse triangle inequality. The strict polygon inequality gives $2\max_i r_i<a+b+S$, hence $L_0<a+b$. Together with the assumption $L_0>|a-b|$, this shows that the argument of $\arccos$ in \eqref{eq:chord} lies in $(-1,1)$. Apply the cosine rule to the two sides $a,b$ and this endpoint distance. The two possible orientations give the stated interval, which also covers reflex angles.
		\end{proof}
		
		\begin{definition}[Forced-angle closure]\label{def:closure}
			For each face, retain the reduced corner expressions $\alpha^*_{uv}=s_{uv}x_{uv}+o_{uv}$ of \eqref{eq:reduced-angles}, where $s_{uv}\in\{-1,1\}$ and $o_{uv}$ is an integral multiple of $60$. Starting with no fixed angle values, repeatedly apply the following rules:
			\begin{enumerate}
				\item[(R1)] In a nondegenerate contracted triangle, the cosine rule determines each remaining corner. If its cosine equals $1/2$, record $s_{uv}x_{uv}+o_{uv}=60$.
				\item[(R2)] If $x_{uv}=60$, record $x_{vu}=180$.
				\item[(R3)] Contract every corner known to satisfy $s_{uv}x_{uv}+o_{uv}=180$, adding the lengths of adjacent sides. Record a contradiction if the resulting polygon has at most two corners or fails the strict polygon inequality.
			\end{enumerate}
			Iterate until no additional fixed angle or contradiction is obtained. A forced value $x_{uv}\notin(0,240)$ contradicts strictness.
		\end{definition}
		
		\begin{lemma}[Soundness of closure]\label{lem:closure}
			Every fixed angle obtained by rules (R1)--(R3) holds in every realization of the marked case. Every recorded polygon contradiction excludes such a realization.
		\end{lemma}
		\begin{proof}
			Induct on the order in which facts are obtained. Initially there are no claims. Suppose all previously recorded angle values hold in a realization. Deleting a known straight corner leaves the same simple polygonal boundary and combines its two collinear sides into one side whose length is their sum. The cosine rule then gives (R1); the contact equation \eqref{eq:contact} gives (R2). Rule (R3) is necessary for any nondegenerate simple polygon by Lemma~\ref{lem:polygon}(1). Each rule therefore preserves necessity. Only finitely many corners can be fixed or removed, so the procedure terminates. A forced value of zero, $240$, or outside that interval contradicts \eqref{eq:contact}.
		\end{proof}

		\section{Finite enumeration}\label{sec:census}
		
		By Lemmas~\ref{lem:minimal}, \ref{lem:faces}, and \ref{lem:three}, a minimal counterexample with fewer than fourteen triangles has a simple, two-connected planar contact graph of even order $t\le12$. Exactly two nonadjacent vertices have degree two; all other vertices have degree three, and distinct three-cycles are vertex-disjoint. We enumerate the graphs and their marked planar embeddings.
		
		The generator first enumerates all connected graphs of maximum degree three on at most twelve vertices, modulo isomorphism. For completeness, observe that every connected graph on at least two vertices has a vertex that is not a cut vertex. Removing it leaves a connected graph of maximum degree three. Conversely, a graph is recovered by adjoining a new vertex to a nonempty subset of at most three vertices of degree at most two. Starting from a single vertex and generating all such extensions is therefore exhaustive. Two implementations use different isomorphism machinery: nauty and the VF2 algorithm in NetworkX.
		
		For each surviving graph, the generator examines all cyclic orders of neighbors and retains only rotation systems of genus zero, as determined by Euler's formula. These give all spherical maps up to graph isomorphism and global reversal of the cyclic orders. Each face is considered as a possible exterior face. For each of the two degree-two vertices, either incident face may contain its free corner. A \emph{marked case} records the rotation system, an outer-face dart, and two pole-face darts. If a map has $F$ faces, all $4F$ choices are retained as \emph{stored records}; \emph{marked classes} identify choices equivalent under graph isomorphism and global reversal.
		
		\begin{table}[htbp]
			\centering
			\caption{Combinatorial census. ``Two-connected'' counts simple two-connected graphs with degree sequence $(2,2,3^{t-2})$ whose degree-two vertices are nonadjacent. ``Admissible'' additionally requires planarity and pairwise vertex-disjoint three-cycles. ``Maps'' counts genus-zero rotation systems modulo graph isomorphism and global reversal. Marked classes use the same symmetries; stored records retain individual markings.}\label{tab:graph-census}
			\begin{tabular}{@{}r r r r r r@{}}
				\toprule
				$t$&Two-connected&Admissible&Maps&Marked classes&Stored records\\
				\midrule
				6&3&3&4&30&64\\
				8&13&9&12&166&240\\
				10&83&47&68&1\,257&1\,632\\
				12&639&272&419&10\,453&11\,732\\
				\bottomrule
			\end{tabular}
		\end{table}
		
		\begin{table}[htbp]
			\centering
			\caption{Exhaustive geometric exclusions for $t<14$. Each stored record is assigned its first successful test. ``Polygon'' counts geometric polygon contradictions; ``Linear certificate'' includes witnesses involving inadmissible forced-angle values.}\label{tab:census}
			\begin{tabular}{@{}r r r r r r@{}}
				\toprule
				$t$&Stored records&Area&Exterior ears&Polygon&Linear certificate\\
				\midrule
				6&64&59&1&4&0\\
				8&240&230&5&4&1\\
				10&1\,632&1\,577&31&18&6\\
				12&11\,732&11\,380&152&175&25\\
				\midrule
				Total&13\,668&13\,246&189&201&32\\
				\bottomrule
			\end{tabular}
		\end{table}
		
		In the census, \emph{Polygon} denotes a degenerate contracted polygon or a violation of its strict polygon inequality. An inadmissible forced-angle value may instead be recorded as a \texttt{FIX} row of a linear certificate. These categories distinguish the first archived witness, not mutually exclusive contradictions.
		
		For $t=4$ the only candidate is $K_4-e$, whose two three-cycles share an edge. Lemma~\ref{lem:three} excludes this case before geometric enumeration. The reductions are summarized in Table~\ref{tab:graph-census}; Table~\ref{tab:census} records the successive geometric tests. All counts can be reproduced from the archived data. Every hypothetical realization is represented by a marked class: its embedding fixes a rotation system, an exterior face, and a face containing each pole. Global reversal may be factored out because reflection preserves strict realizability.

		\section{Exact exclusion certificates}\label{sec:certs}
		
		For a stored record, let $k$ be the number of unit sides of the outer boundary, $p$ the number of exterior poles, and $f_{\rm odd}$ the number of bounded odd-sided white faces. The first two computational tests are direct consequences of Lemmas~\ref{lem:area} and \ref{lem:ear}:
		\begin{align}
			(t+f_{\rm odd})\sqrt3&\le4M(k),\label{eq:test-area}\\
			(t-p+f_{\rm odd})\sqrt3&\le4M(k-p).\label{eq:test-ear}
		\end{align}
		The implementation rejects a case only when rational strict comparisons establish the failure of \eqref{eq:test-area} or \eqref{eq:test-ear}. Equality is not a rejection; in particular $4M(6)=6\sqrt3$. The rational estimates
		\begin{equation}\label{eq:rat-bounds}
			\frac{173}{100}<\sqrt3,\qquad
			\frac{173205}{100000}<\sqrt3,\qquad
			\frac{\sqrt3}{4}<\frac{433013}{10^6}
		\end{equation}
		are sufficient. Upper bounds for $M(k)$ follow from $\pi>3.14159$ and $\tan x\ge x+x^3/3$ for $0\le x<\pi/2$, together with simple exact bounds for small $k$. Thus no floating-point calculation determines whether a case is rejected.
		
		Table~\ref{tab:census} shows that these two filters discard 13\,435 cases, leaving 233. For 201 of these cases, the closure rules of Definition~\ref{def:closure} yield a polygon with fewer than three corners or with side lengths violating the strict polygon inequality. Each of these \emph{semantic certificates} identifies a face and a verifiable contradiction.
		
		The other 32 cases are excluded by rational linear certificates in the contact angles $x_{uv}$. For each fixed $t$, record indices are zero-based positions in \path{data/cases_tXX.json}, with matching certificates at the same positions in \path{data/certs_tXX.json}, where \texttt{XX} is the two-digit triangle count. The tests run in the order area, exterior ears, contracted polygon, linear certificate. A \texttt{CDEG} certificate stores a witness dart for a degenerate contracted polygon. Each linear row expresses a necessary geometric condition:
		\begin{itemize}
			\item \texttt{CONTACT} and \texttt{FACE}: equations \eqref{eq:contact} and \eqref{eq:face-angle};
			\item \texttt{FIX}: a forced rational angle obtained by rules (R1)--(R2), including values at or outside the admissible interval (which already imply a contradiction);
			\item \texttt{QEQ}: equality of opposite base angles in a kite or rhombus;
			\item \texttt{HI}: a rational bound from the triangle cosine rule;
			\item \texttt{CHL}: a rational bound derived from \eqref{eq:chord}.
		\end{itemize}
		
		\begin{lemma}[Exact linear exclusion]\label{lem:motzkin}
			Suppose the geometric rules yield equalities $a_i\cdot x=b_i$, inequalities $g_j\cdot x\ge c_j$, and $0<x_h<240$ for every variable. Choose rational $\lambda_i$, arbitrary in sign, and rational $\rho_j\ge0$. Set
			\begin{align*}
				w&=\sum_i\lambda_i a_i-\sum_j\rho_jg_j,\\
				B&=\sum_i\lambda_i b_i-\sum_j\rho_jc_j
				+240\sum_h\max(-w_h,0).
			\end{align*}
			If $B<0$, or if $B\le0$ and $w\ne0$, then the marked case is not geometrically realizable.
		\end{lemma}
		\begin{proof}
			For a feasible $x$, direct substitution gives \[ \sum_{w_h>0} w_hx_h +\sum_{w_h<0}(-w_h)(240-x_h) =B-\sum_j\rho_j(g_j\cdot x-c_j)\le B. \] The left side is nonnegative; it is strictly positive if $w\ne0$, by the strict angular bounds. Either hypothesis is a contradiction.
		\end{proof}
		
		All certificate multipliers are rational, and all rows are independently reconstructed from the marked map before the above inequality is checked. Eight certificates use $\arccos(1/4)<90^\circ$, which follows immediately from $1/4>0$. The only other trigonometric constant occurring in a retained angular certificate is $\arccos(1/2)=60^\circ$. Hence the angular certificate checker needs no approximation of $\pi$ or transcendental function.
		
		\begin{example}[A collapsed polygon]
			The stored record of index 36 for $t=6$ has the bounded face with dart cycle \[ (2,4),\ (4,5),\ (5,3),\ (3,2). \] The exact closure forces $x_{45}=x_{32}=180^\circ$. Removing these two straight corners leaves only the two contact corners $(2,4)$ and $(5,3)$, with successive side lengths $(2,2)$. No nondegenerate simple polygon can have two corners, so the record is excluded by a \texttt{CDEG} certificate whose witness is the dart $(2,4)$.
		\end{example}
		
		\Needspace{16\baselineskip}
		\begin{example}[An exact linear contradiction]
			For $t=8$, the stored record of index 146 can be excluded using nineteen linear equations: eleven contact-angle equations, two face-angle equations, and six equalities between opposite angles in unit rhombi. Using the vertex labels of the archived marked map, their row identifiers, in order, are \[
			\begin{aligned}
				\texttt{CONTACT}:\quad
				&(04,07,15,16,23,24,27,35,36,45,67),\\
				\texttt{FACE}:\quad &(07,16),\\
				\texttt{QEQ}:\quad
				&((04,27),(42,70),(23,67),\\
				&\hspace{1.1em}(36,72),(24,53),(45,32)).
			\end{aligned}
			\] Here $uv$ denotes the ordered dart $(u,v)$, and $(uv,wz)$ identifies the equation $x_{uv}=x_{wz}$. The two face rows have the following dart cycles and parameters: \[
			\begin{aligned}
				(07):\quad &(07,76,61,15,54,40),&\ell=6,\quad p_W=2,\quad&\text{exterior};\\
				(16):\quad &(16,63,35,51),&\ell=4,\quad p_W=0,\quad&\text{bounded}.
			\end{aligned}
			\] Their face-equation right-hand sides are $1200^\circ$ and $360^\circ$, respectively, by \eqref{eq:face-angle}. Multiply the three groups of equations respectively by \[
			\begin{aligned}
				\lambda_C&=(-1,-1,-1,-1,1,1,1,-1,-1,-1,-1),\\
				\lambda_F&=(1,1),\\
				\lambda_Q&=(1,-1,-1,1,-1,1).
			\end{aligned}
			\] All angle coefficients cancel, while the constants sum to $360^\circ$. Thus the equations imply $0=360^\circ$, a contradiction. This certificate requires no angular inequalities. The electronic record contains the negatives of the displayed multipliers, so that its verifier obtains $B=-360$ in Lemma~\ref{lem:motzkin}.
		\end{example}
		
		\begin{proposition}[Finite exclusion]\label{prop:finite}
			No strict configuration having the structural properties (1)--(3) of Lemma~\ref{lem:minimal} has $t\le12$ triangles.
		\end{proposition}
		\begin{proof}
			The exhaustive graph and embedding generation in Section~\ref{sec:census} produces the marked cases in Table~\ref{tab:census}. By the first two verified tests, 13\,435 cases violate a necessary area inequality. Of the 233 remaining cases, the supplementary standard-library verifier validates 201 contracted-polygon contradictions and 32 rational certificates, using Lemmas~\ref{lem:angles}, \ref{lem:three}, \ref{lem:polygon}, and \ref{lem:motzkin}. Each marked case is therefore excluded. The enumeration is exhaustive by the vertex-extension argument and the complete enumeration of planar rotation systems, outer faces, and pole sectors.
		\end{proof}
		
		\begin{remark}[Reproducibility]\label{rem:reprod}
			The electronic supplement~\cite{WinklerSupplement2026} provides four case lists, four certificate lists, and the standard-library verifier \path{check233_stdlib.py} for Proposition~\ref{prop:finite}. The verifier checks vertex labels, degree sequences, adjacency, two-connectivity, vertex-disjointness of three-cycles, planar facial cycles, graph6 encodings, pole incidence, certificate indices, and record counts before applying the geometric rules. Independent graph generators, one without nauty, and bidirectional comparisons of planar maps and marked classes establish the completeness of the lists. A separate implementation independently checks all 233 geometric certificates; its programs and logs are supplied. The exact constructions of the fourteen-triangle configurations are verified by \path{check_D_stdlib.py}, also using the Python standard library.
		\end{remark}

		\section{Exact fourteen-triangle constructions}\label{sec:construction}
		
		We give exact coordinates for the double kite and the reverse double kite. Let
		\begin{align*}
			s&=(0,0),&
			a_1&=\left(-\frac12,\frac{\sqrt{15}}2\right),&
			a_2&=\left(\frac12,\frac{\sqrt{15}}2\right),\\
			p&=\left(0,\frac{\sqrt{15}+\sqrt3}{2}\right),&
			b_1&=R_{60^\circ}(a_1),&
			b_2&=R_{-60^\circ}(a_2),
		\end{align*}
		where $R_\varphi$ denotes rotation through $\varphi$ about the origin. Both $sa_1b_1$ and $sa_2b_2$ are equilateral triangles of side length two. For any such large triangle $ABC$, let $\mathcal{C}(A,B,C)$ denote the three unit corner triangles with vertices \[ \left(A,\frac{A+B}{2},\frac{A+C}{2}\right),\quad \left(B,\frac{B+A}{2},\frac{B+C}{2}\right),\quad \left(C,\frac{C+A}{2},\frac{C+B}{2}\right). \] The seven-triangle \emph{kite} is
		\begin{equation}\label{eq:kite}
			K=\mathcal{C}(s,a_1,b_1)\;\cup\;\mathcal{C}(s,a_2,b_2)\;\cup\;\{\triangle a_1a_2p\}.
		\end{equation}
		Here the unions are unions of families of triangles. Let $\sigma$ be reflection in the line through $p$ and $b_2$. Set
		\begin{equation}\label{eq:dk}
			D=K\cup\sigma(K).
		\end{equation}
		The \emph{reverse double kite} is obtained by the half-turn $\tau(x)=p+b_2-x$, whose center is the midpoint of $p$ and $b_2$:
		\begin{equation}\label{eq:rdk}
			D_{\mathrm{rev}}=K\cup\tau(K).
		\end{equation}
		Both configurations are shown in Figure~\ref{fig:kite-pair}.
		
		\begin{figure}[htbp]
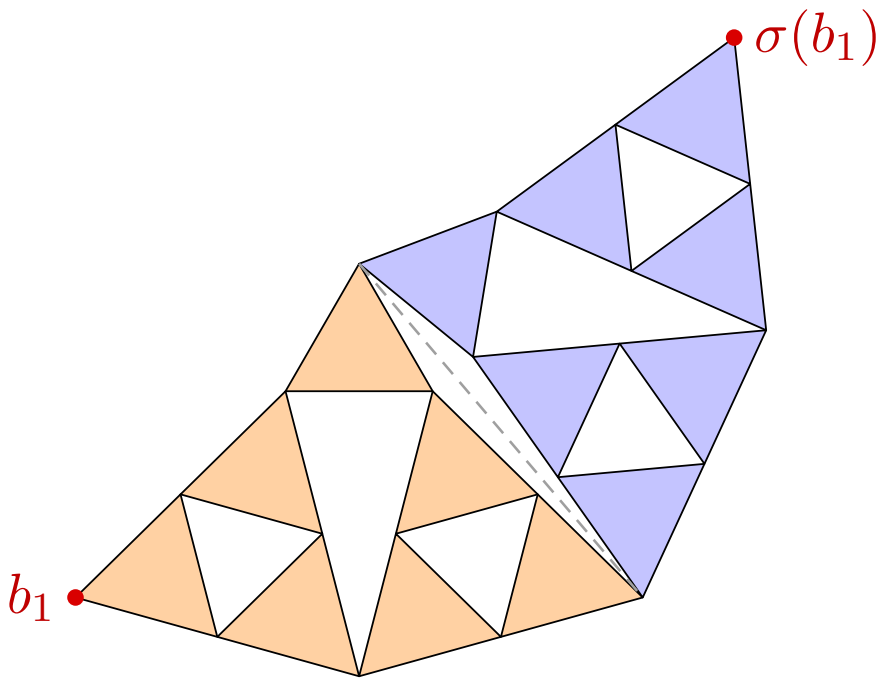

			\centering
			\resizebox{.4\textwidth}{!}{\DKDoubleKite}
			\caption{The double kite $D$, obtained by reflecting a seven-triangle kite across the dashed line (Proposition~\ref{prop:existence}). Orange and blue distinguish the two constituent kites, as in Figure~\ref{fig:kite-pair}. The red points are the poles $b_1$ and $\sigma(b_1)$. The drawing uses decimal evaluations of the exact coordinates in \eqref{eq:kite}--\eqref{eq:dk}.}\label{fig:dk}
		\end{figure}
		
		\begin{proposition}\label{prop:existence}
			Each of the families $D$ and $D_{\mathrm{rev}}$ defined by \eqref{eq:kite}--\eqref{eq:rdk} is a strict connected configuration with fourteen unit triangles and exactly two poles.
		\end{proposition}
		\begin{proof}
			The two large triangles in \eqref{eq:kite} have side length two, hence each of their six corner triangles has unit sides. Moreover, $a_1a_2p$ is equilateral with unit sides. Reflection and the half-turn preserve these lengths. Thus both families consist of fourteen unit triangles.
			
			All coordinates lie in $\mathbb{Q}(\sqrt3,\sqrt5)$. In either family, exact coordinate comparisons give 22 distinct vertices, 20 of multiplicity two and precisely two of multiplicity one: $b_1,\sigma(b_1)$ for $D$, and $b_1,\tau(b_1)$ for $D_{\mathrm{rev}}$. Orientation and segment-intersection predicates verify that any two distinct triangles are disjoint or meet at a common vertex, and that the 20 contacts give a connected contact graph. The supplementary standard-library program \path{check_D_stdlib.py} checks these claims for both configurations using exact fractions in $\mathbb{Q}(\sqrt3,\sqrt5)$. Orientation signs are determined by rational isolating intervals, without numerical tolerances.
		\end{proof}
		
		\begin{proof}[Proof of Theorem~\ref{thm:main}]
			Suppose a strict configuration with one or two poles and fewer than fourteen triangles exists. Then the class of strict configurations with one or two poles and at most thirteen triangles is nonempty. Choose a minimum member. By Lemma~\ref{lem:minimal}, it has even triangle count at most twelve and satisfies all structural assumptions of Proposition~\ref{prop:finite}, a contradiction. Proposition~\ref{prop:existence} shows that both double-kite constructions attain the lower bound.
		\end{proof}
		
		\begin{corollary}\label{cor:one-pole}
			Every strict one-pole configuration contains at least fifteen triangles, if such a configuration exists.
		\end{corollary}
		\begin{proof}
			By Theorem~\ref{thm:main} a one-pole configuration contains at least fourteen triangles. Its number of triangles is odd by Lemma~\ref{lem:count}.
		\end{proof}
		
		\begin{corollary}[Small edge cuts]\label{cor:edgecut}
			Let $\T$ be a strict configuration without poles, and let $H$ be its contact graph. If $H$ is not three-edge-connected, then $\T$ contains at least twenty-eight unit triangles. If $H$ has a bridge, then $\T$ contains at least thirty triangles.
		\end{corollary}
		\begin{proof}
			The contact graph is cubic. A bridge separates it into two connected subgraphs; each determines a strict one-pole configuration after the contact is cut. Corollary~\ref{cor:one-pole} gives at least fifteen triangles on either side. Otherwise, a minimal two-edge cut separates $H$ into two connected subgraphs, each defining a strict two-pole configuration. The two newly free contacts are distinct because the original configuration is strict. Theorem~\ref{thm:main} gives at least fourteen triangles on either side.
		\end{proof}
		
		\begin{corollary}[Decompositions into two-pole configurations]\label{cor:ring}
			Suppose a strict four-regular unit-triangle configuration is the union of $k\ge 2$ strict two-pole configurations whose triangle families are pairwise disjoint. Then it contains at least $14k$ unit triangles. In particular, a decomposition into three such components requires at least forty-two triangles, and this bound is attained by both configurations in Figure~\ref{fig:fortytwo}.
		\end{corollary}
		\begin{proof}
			The number of unit triangles is the sum of the sizes of the $k$ component families. Each has at least fourteen triangles by Theorem~\ref{thm:main}, which proves the lower bound.
			
			To verify sharpness for $k=3$, first consider the double kite $D$ from \eqref{eq:dk}, with poles $b_1$ and $q=\sigma(b_1)$. Let $R$ be counterclockwise rotation through $120^\circ$ about the origin, and define \[ c=(I-R)^{-1}(q-Rb_1),\qquad \rho(x)=c+R(x-c). \] Then $\rho(b_1)=q$. The three families $D$, $\rho(D)$, and $\rho^2(D)$ form Harborth's ring. The exact computation in \path{src/check_ring42_DDD.py} establishes strictness, connectivity, and the fact that each of its sixty-three vertices belongs to exactly two triangles. The clockwise rotation instead makes the three copies overlap.
			
			For the asymmetric example, take two copies of $D$ and one copy of $D_{\mathrm{rev}}$ from \eqref{eq:rdk}. Their squared pole separations, calculated from the exact coordinates, are \[ d_D^2=\lvert b_1-\sigma(b_1)\rvert^2=\frac{687+279\sqrt5}{38},\qquad d_{\mathrm{rev}}^2=\lvert b_1-\tau(b_1)\rvert^2=\frac{57+9\sqrt5}{2}. \] In particular, $0<d_{\mathrm{rev}}<2d_D$. Put \[ J_0=(0,0),\qquad J_1=(d_D,0),\qquad J_2=\left(\frac{d_{\mathrm{rev}}^2}{2d_D}, \sqrt{d_{\mathrm{rev}}^2-\frac{d_{\mathrm{rev}}^4}{4d_D^2}}\right). \] The sides $J_0J_1$, $J_1J_2$, and $J_2J_0$ have lengths $d_D,d_D,d_{\mathrm{rev}}$. Take the three ordered components $(D,D,D_{\mathrm{rev}})$, with ordered poles $(b_1,\sigma(b_1))$, $(b_1,\sigma(b_1))$, and $(b_1,\tau(b_1))$, respectively. Map each component by the unique orientation-preserving isometry sending its first pole to $J_i$ and its second to $J_{i+1}$ (indices modulo three), without reflecting the component. This is exactly the placement with \texttt{mirror=(0,0,0)} and \texttt{orient=(0,0,0)} in \path{src/check_ring42_balls.py}. The interval certificates establish that these three images are a strict ring of forty-two triangles and sixty-three doubly occupied vertices. Hence both rings attain the bound.
		\end{proof}

		\section{Verification and scope}\label{sec:scope}
		Theorem~\ref{thm:main} gives the sharp minimum for strict two-pole configurations, but does not classify the configurations attaining it. Shared edges or contacts in the relative interior of a triangle side are excluded by Definition~\ref{def:model}. Corollary~\ref{cor:edgecut} applies to strict four-regular configurations whose contact graphs have an edge cut of size one or two; Corollary~\ref{cor:ring} gives a componentwise bound for specified decompositions. Whether two two-pole components can form a strict four-regular configuration remains open, as does the minimum size among configurations with three-edge-connected contact graphs.
		
		The file \path{DoubleKite_Minimality_Supplement_v1.0_2026-11-09.zip} is the electronic supplement archived at Zenodo (version 1.0, November 9, 2026; \url{https://doi.org/10.5281/zenodo.24093082}). It contains every marked case, the exclusion certificates and their verifier, the independent enumeration and audit programs, and tests of the exact fourteen-triangle constructions and both forty-two-triangle rings. The DDD ring is verified using exact field arithmetic; the asymmetric DDR ring is verified by certified interval arithmetic. The accompanying README specifies dependencies and reproduction commands.
		
		\paragraph{Computational and AI assistance.}
		Claude (Anthropic, Opus 5.5), ChatGPT (OpenAI, GPT-6), and Gemini (Google, Argon) were used during computational exploration, certificate verification, and preparation of the manuscript. Claude contributed to graph enumeration and certificate construction and later performed a code-based completeness audit using a separate nauty-free generator and independently implemented geometric tests. ChatGPT assisted with the area-filter and certificate-checker audits, exact-coordinate verification, and exposition. Gemini provided a referee-style critique and a partial independent checking program. The graph-reduction strategy, geometric exclusion rules, and certificate representation were developed in an iterative process with Claude and ChatGPT under the author's direction. The author checked and takes responsibility for the mathematical claims, the reproducible programs, and the submitted text.

	\end{document}